\documentclass[10pt]{amsart}
\usepackage{graphics,graphicx, xcolor}
\usepackage{caption}
\usepackage{amsfonts}
\usepackage{epsfig}
\usepackage{verbatim}
\usepackage{amssymb}
\usepackage{latexsym}
\usepackage{amsmath, amsthm, amssymb}

\newtheorem{corollary}{Corollary}[section]

\newtheorem{lemma}[corollary]{Lemma}

\newtheorem{theorem}[corollary]{Theorem}
\newcommand{\mylabel}[1]{\label{#1}
            \ifx\undefined\stillediting
            \else \fbox{$#1$}\fi }
\newcommand{\BE}{\begin{equation}}

\newcommand{\EEQ}{\end{equation}}
\newcommand{\rfb}[1]{\mbox{\rm
   (\ref{#1})}\ifx\undefined\stillediting\else:\fbox{$#1$}\fi}

\newfont{\Blackboard}{msbm10 scaled 1200}

\newfont{\roma}{cmr10 scaled 1200}

\def\CC{\rm \hbox{C\kern-.56em\raise.4ex
         \hbox{$\scriptscriptstyle |$}\kern+0.5 em }}

\makeatletter

\@addtoreset{equation}{section}
\makeatother

\newcommand{\mm}    {{\hbox{\hskip 0.5pt}}}

\newcommand{\bluff} {{\hbox{\raise 15pt \hbox{\mm}}}}
\makeatletter
\def\section{\@startsection {section}{1}{\z@}{-3.5ex plus -1ex minus
    -.2ex}{2.3ex plus .2ex}{\large\bf}}
\makeatother
\def\be{\begin{equation}}
\def\ee{\end{equation}}

\input{epsf}

\begin{document}
\thispagestyle{empty}
\title[Thermoelastic system with delay]{A note on exponential stability of a thermoelastic system with internal delay}
\author{Smain Moulay Khatir}
\address{LR Analysis and Control of PDEs, LR 22ES03, Department of Mathematics, Faculty of Sciences of Monastir, University of Monastir, Tunisia, and
Laboratory of Analysis and Control of Partial Differential Equations, Djillaly Liabes University, Sidi Bel Abbes, Algeria.}
\email{s.moulay\_khatir@yahoo.fr}

\author{Farhat Shel}
\address{LR Analysis and Control of PDEs, LR 22ES03, Department of Mathematics, Faculty of Sciences of Monastir, University of Monastir, Tunisia}
\email{farhat.shel@fsm.rnu.tn} 

\begin{abstract} The presence of a delay in a thermoelastic system destroys the well-posedness and the stabilizing effect of the heat conduction. To avoid this problem we add to the system, at the delayed equation, a Kelvin-Voigt damping. In this note we point on the exponential stability of such system in order to improve the mean  result in our paper Well-posedness and exponential stability of a thermoelastic system with internal delay (Applicable Analysis J 101, 4851-4865, 2022). We use a frequency domain method in the proof of stability.
\end{abstract}


\maketitle

\tableofcontents

\section{Introduction}

Let $\Omega:=(0,\ell)$, $\ell>0$. In \cite{KhSh21}, we considered the following thermoelastic system with delay
\begin{equation}
\left\{ 
\begin{tabular}{l}
$u_{tt}(x,t)-\alpha u_{xx}(x,t-\tau)-\beta u_{xxt}(x,t)+\gamma \theta _x(x,t)=0,\;\;\text{ in } \Omega\times (0,\infty ),$ \\
$\theta_{t}(x,t)-\kappa\theta_{xx}(x,t)+\gamma
u_{xt}(x,t)=0,\;\;\;\;\;\;\;\;\;\;\;\;\;\;\;\;\;\;\;\;\;\;\;\;\;\;\;\;\;\;\text{ in } \Omega \times (0,\infty ),$\\
$u(0,t)=u(\ell ,t)=0,\;\;\;\;\;\;\;\;\;\;\;\;\;\;\;\;\;\;\;\;\;\;\;\;\;\;\;\;\;\;\;\;\;\;\;\;\;\;\;\;\;\;\;\;\;\;\;\;\;\;\;\;\;\;\;\;\text{ in } (0,\infty ),$ \\ 
$\theta_x(0,t)=\theta_x(\ell ,t)=0,\;\;\;\;\;\;\;\;\;\;\;\;\;\;\;\;\;\;\;\;\;\;\;\;\;\;\;\;\;\;\;\;\;\;\;\;\;\;\;\;\;\;\;\;\;\;\;\;\;\;\;\;\;\text{ in }(0,\infty ),$ \\ 
$u_{x}(x,t-\tau)=f_0(x,t-\tau),\;\;\;\;\;\;\;\;\;\;\;\;\;\;\;\;\;\;\;\;\;\;\;\;\;\;\;\;\;\;\;\;\;\;\;\;\;\;\;\;\;\;\;\;\;\;\text{ in } \Omega\times (0,\tau ),$ \\
$u(x,0)=u_0(x), u_t(x,0)=u_1(x), \theta (x,0)=\theta _0(x),\;\;\;\;\;\;\;\;\;\;\;\;\text{ in } \Omega$
\end{tabular}
\right.  \label{d1.0}
\end{equation}
where $\alpha , \gamma , \kappa , \beta $ are some positive constants. The functions $u=u(x,t)$ and $\theta =\theta (x,t)$ describe respectively the displacement and the temperature difference, with $x \in \Omega$ and $t \geq 0.$ Moreover,
 $\tau >0$ is the time delay. The initial data $(u_0,u_1,  f_0, \theta_0)$ belongs to a suitable space.

Recall that the term  $-\beta u_{xxt}(x,t)$ serves to avoid the problem of illposedness and instability caused by the delay present at $u_{xx}(x,t-\tau)$. In this note we try to improve the result of exponential stability given in \cite{KhSh21}, Theorem 3.1,  by specifying a concrete lower bound of the damping coefficient $\beta$ in terms of $\alpha\tau$, as in \cite{Amm15'}. The method used a result due to Gearhard, Prûss and Huang which characterize the exponential stability of a semigroup by means of the frequency domain approach.

 We introduce, as in \cite{Amm15'} and \cite{KhSh21} , the new variable 
 \begin{equation*}
 z(x,\rho ,t)=u_x(x,t-\tau \rho),\;\;\;\;\;\;\text{ in } \Omega\times (0,1)\times (0,\infty ),
 \end{equation*}
 The new function $z(x,\rho,t)$ satisfies
 \begin{eqnarray*}
 \tau z_t(x,\rho,t)+z_{\rho}(x,\rho,t)&=&0,\;\;\;\;\;\;\;\;\;\;\;\;x\in \Omega,\;\; \rho \in (0,1),\;\; t \in (0,+\infty),\\
 z(x,0,t)&=&u_x(x,t),\;\;\;x\in  \Omega,\;\; t \in (0,+\infty).
 \end{eqnarray*}
 
Then, problem (\ref{d1.0}) takes the form
\begin{eqnarray}
u_{tt}(x,t)-\alpha z_{x}(x,1,t)-\beta u_{xxt}(x,t)+\gamma \theta _x(x,t)=0,&&\text{ in } \Omega\times (0,\infty ), \label{e1}\\
\tau z_t(x,\rho,t)+z_{\rho}(x,\rho,t)=0,&&\text{ in } \Omega\times(0,1)\times (0,+\infty), \label{e2}\\
\theta_{t}(x,t)-\kappa\theta_{xx}(x,t)+\gamma
u_{xt}(x,t)=0,&&\text{ in } \Omega\times (0,\infty ), \label{e3}\\
u(0,t)=u(\ell ,t)=0,&&\text{ in } (0,\infty ), \\ 
\theta_x(0,t)=\theta_x(\ell ,t)=0,&&\text{ in }(0,\infty ), \label{e4}
\end{eqnarray}
\begin{eqnarray} 
z(x,0,t)=u_{x}(x,t),&&\text{ in } \Omega\times (0,\infty ), \label{e5}\\
u(x,0)=u_0(x), u_t(x,0)=u_1(x), \theta (x,0)=\theta _0(x),&&\text{ in } \Omega, \label{e6}\\
 z(x,\rho, 0)=f_0(x,-\tau \rho),&&\text{ in } \Omega\times (0,1 ). \label{e7}
\end{eqnarray}

Without loss of generality, we assume that $\int_\Omega \theta(x,t)dx=0$. Otherwise, we can make the substitution $\tilde{\theta}(x,t)=\theta (x,t)-\frac{1}{\ell}\int_\Omega \theta_0(x)dx,$ in fact $(u,v,z,\theta)$ and $(u,v,z,\tilde{\theta})$ satisfy the same system (\ref{e1})-(\ref{e7}).

Let 
\begin{equation*}
\mathcal{H}=\left\{(f,g,p,h)\in H^{1}_0(\Omega)\times L^2(\Omega)\times L^2(\Omega \times (0,1)) \times L^2(\Omega)\mid \int_\Omega h(x)dx=0 \right\}.
\end{equation*}
Equipped with the following inner product: for any $U_k=(f_k,g_k,p_k,h_k) \in \mathcal{H},\;\;k=1,2,$
\begin{equation*}
\left\langle U_1,U_2\right\rangle _{\mathcal{H}}=\int_\Omega \left(\alpha f_{1x}(x)f_{2x}(x)+ g_1(x)g_2(x)+h_1(x)h_2(x) \right)dx+\xi\int_\Omega \int_0^1p_1(x,\rho)p_2(x,\rho)d\rho dx, 
\end{equation*}
$\mathcal{H}$ is a Hilbert space.

Define
\begin{equation*}
U:=(u,u_t,z,\theta)
\end{equation*} 
then, problem (\ref{d1.0}) can be formulated as a first order system of the form
\begin{equation}
\left\{ 
\begin{array}{c}
U^\prime =\mathcal{A}U, \\ 
U(0)=(u_0,u_1,f_0(.,-.\tau), \theta_0)
\end{array}
\right.  \label{d500}
\end{equation}
where the operator $\mathcal{A}$ is defined by
\begin{equation*}
\mathcal{A}\left( 
\begin{array}{c}
u \\ 
v \\ 
z \\ 
\theta
\end{array}
\right) =\left( 
\begin{array}{c}
v \\ 
(\alpha z(.,1)+\beta v_{x})_x-\gamma \theta _x \\ 
-\frac{1}{\tau}z_\rho \\ 
-\gamma v_{x}+\kappa \theta _{xx} 
\end{array}
\right)
\end{equation*}
with domain
\begin{equation*}
\mathcal{D}(\mathcal{A})=\left\{ 
\begin{array}{c}
U=(u,v,z,\theta)\in \mathcal{H}\cap \left[ H^{1}_0(\Omega)\times H^{1}_0(\Omega)\times L^2(\Omega ; H^1(0,1))\times  H^2(\Omega)\right]   \mid \\ 
\;\;z(.,0)=u_x,\;\;\;\theta_x(0)=\theta_x(\ell)=0\;\;\;\text{and} \;\; 
(\alpha z(.,1)+\beta v_{x}) \in H^1(\Omega)
\end{array}
\right\}
\end{equation*}
in the Hilbert space $\mathcal{H}.$

It can be proved (see proof of Lemma 2.1 in \cite{KhSh21})
\begin{lemma}
If $\xi=\frac{2\tau \alpha^2}{\beta}$, then  $\mathcal{A}-mId$, with $m=2\frac{\alpha^2}{\beta}$, is dissipative maximal and satisfies
\begin{equation}
Re \left( \left\langle \mathcal{A}U,U\right\rangle _{\mathcal{H}}\right)\leq  -\frac{1}{2}\beta \int_\Omega |v_{x}(x)|^2dx+m \int_\Omega  |u_x(x)|^2dx-\kappa\int_\Omega |\theta _{x}(x)|^2dx.\label{ccsp}
\end{equation}
Moreover, the operator $\mathcal{A}$ generates a $\mathcal{C}_0$-semigroup on $\mathcal{H}$.
 \end{lemma}
  
\section{Exponential stability}
Now, we prove that under the condition $\beta \geq \alpha \tau$, the semigroup $e^{t\mathcal{A}}$ is exponentially stable.
\begin{theorem} \label{cc4}
If $\xi=2\frac{\tau \alpha^2}{\beta}$ and $ \alpha\tau \leq \beta$ then there exists two positive constants $C$ and $w$ such that the semigroup $e^{t\mathcal{A}}$ satisfies the following estimate:
$$\|e^{t\mathcal{A}}\|_{\mathcal{L}(\mathcal{H})}\leq Ce^{-wt},\;\;\;\;\;\forall\,t>0.$$ 
\end{theorem}
\begin{proof}
We will use a frequency domain method, based on the following result due to Gearhart, Prüss and Huang \cite{Liu99}.
\begin{lemma}\label{cc3}
A $\mathcal{C}_0$-semigroup $e^{t\mathcal{B}}$ on a Hilbert space $H$ satisfies
\begin{equation*}
\|e^{t\mathcal{B}}\|_{\mathcal{L}(H)}\leq Ce^{-wt},\;\;\;\;\;\forall\,t>0,
\end{equation*}
for some positive constants $C$ and $w$, if and only if,
\begin{equation*}
Re(\lambda)\leq 0,\;\;;\forall\,\lambda \in \sigma(\mathcal{B})
\end{equation*}
and
\begin{equation*}
\sup_{Re(\lambda)>0}\|(\lambda I-\mathcal{B})^{-1}\|_{\mathcal{L}(H)}<\infty,
\end{equation*}
where $\sigma(\mathcal{B})$ denotes the spectrum of the operator $\mathcal{B}$.
\end{lemma}
\begin{lemma}\label{cc1}
For $\xi=2\frac{\tau \alpha^2}{\beta}$ and $\alpha\tau\leq\beta$,  it holds
\begin{equation}
Re(\lambda)\leq 0,\;\;\forall\,\lambda \in \sigma(\mathcal{A}). \label{c1}
\end{equation}
\end{lemma}
\begin{proof}
Let $\lambda=a+\mathbf{i}b \in \mathcal{C}_0:=\left\lbrace \lambda \in \mathbb{C}\mid  Re(\lambda)>0 \right\rbrace $, with $a=Re(\lambda)$ and $b=Im(\lambda)$ and $F=(f,g,p,h)\in \mathcal{H}$. We look for $U=(u,v,z,\theta)$ such that
\begin{equation}
(\lambda I-\mathcal{A})U=F, \label{sp111}
\end{equation}
i.e.
\begin{equation}
\left\{ 
\begin{tabular}{l}
$\lambda u- v=f,\;\;\text{in}\;\;H_0^1(\Omega),$ \\
$\lambda v-\left(\alpha z(.,1)+\beta v_{x} \right)_x+\gamma \theta_x =g,\;\;\text{in}\;\;L^2(\Omega) ,$\\
$\lambda z+ \frac{1}{\tau}z_{\rho} =p,\;\;\text{in}\;\;L^2(\Omega \times (0,1)),$ \\ 
$\lambda \theta + \gamma v_{x} -\kappa \theta_{xx}= h,\;\;\text{in}\;\;L^2(\Omega).$
\end{tabular}
\right.  \label{max1sp1}
\end{equation}
Again,
\begin{equation*}
z(x,\rho)=e^{-\lambda \tau \rho}u_{x}(x)+ \tau e^{-\lambda \tau \rho}\int_0^\rho p(s)e^{\lambda \tau s}ds,\label{max222}
\end{equation*}
and, in particular
\begin{equation}
z(x,1)=e^{-\lambda \tau}u_{x}+ z_{0}, \label{max33}
\end{equation}
with $z_{0} \in L^2(\Omega)$ defined by
\begin{equation*}
z_{0}=\tau e^{-\lambda \tau}\int_0^1 p(s)e^{\lambda \tau s}ds. \label{max44}
\end{equation*}
Multiplying (\ref{max1sp1})$_2$ and (\ref{max1sp1})$_4$ respectively by $\lambda w\in H^1_0(\Omega)$ and $\varphi\in H^1(\Omega)$, with $\varphi_x(0)=\varphi_x(\ell)=0$, we obtain after some integrations by parts taking into account (\ref{max1sp1})$_1$, (\ref{max33}),  and boundary conditions,
\begin{equation}
\overline{\lambda}\lambda^2 \int_\Omega u\overline{w}
 dx+\left(\alpha \overline{\lambda}e^{-\lambda \tau}+|\lambda|^2 \beta \right)  \int_\Omega u_x\overline{w_x}
 dx \gamma \overline{\lambda}\int_\Omega \theta_{x}  \overline{w} dx =\overline{\lambda} \int_\Omega (g+\lambda f)\overline{w}
 dx +\overline{\lambda} \int_\Omega (\beta f_{x}-\alpha z_{0})\overline{w_{x}} dx\label{eq2.6}
\end{equation}
and
\begin{equation}
\lambda\int_\Omega \theta\overline{\varphi}
 dx + \kappa  \int_\Omega \theta_x\overline{\varphi_x}
 dx -\gamma \overline{\lambda}\int_\Omega \theta_{x} {\color{red}} \overline{w} dx=\int_\Omega (h-\gamma f) \overline{\theta} dx. \label{eq2.7} 
\end{equation}
Summing (\ref{eq2.6}), and (\ref{eq2.7}),  we get
\begin{equation}
B\left((u,\theta),(w,\varphi) \right)=\Phi(w,\varphi) \label{max99}
\end{equation}
with
\begin{eqnarray*}
B\left( (u,\theta),(w,\varphi)\right)&=& \overline{\lambda}\lambda^2 \int_\Omega u\overline{w}
 dx+\left(\alpha \overline{\lambda}e^{-\lambda \tau}+|\lambda|^2 \beta \right)  \int_\Omega u_x\overline{w_x}
 dx + \lambda\int_\Omega \theta\overline{\varphi}
 dx \notag \\&+& \kappa  \int_\Omega \theta_x\overline{\varphi_x}
 dx +\mathbf{i} Im\left(\gamma \overline{\lambda}\int_\Omega \theta_{x}  \overline{w} dx  \right) \label{max7sp111}
\end{eqnarray*}
and
\begin{equation*}
\Phi((w,\varphi))=\overline{\lambda} \int_\Omega (g+\lambda f)\overline{w}
 dx +\overline{\lambda} \int_\Omega (\beta f_{x}-\alpha z_{0})\overline{w_{x}} dx+\int_\Omega (h-\gamma f) \overline{\theta} dx 
\end{equation*}
The space $$\mathcal{F}:=\left\lbrace (w,\varphi) \in H^1_0(\Omega)\times H^1(\Omega)\mid \varphi_x(0)=\varphi_x(\ell)=0 \right\rbrace,$$
equipped with the inner product
$$\langle (w_1,\varphi_1),(w_2,\varphi_2) \rangle _\mathcal{F}=\int_\Omega \left(w_1\overline{w_2}+w_{1x}\overline{w_{2x}}+\varphi_1\overline{\varphi_2}+\varphi_{1x}\overline{\varphi_{2x}} \right)dx, $$
is a Hilbert space; the sesquilinear form $B$ on $\mathcal{F}\times \mathcal{F}$ and the antilinear form $\Phi$ on $\mathcal{F}$ are continuous. Now, we will study carefully the coercivity of $B$. Let $(w,\varphi) \in \mathcal{F},$ one has
$$Re\left( B\left( (w,\varphi),(w,\varphi)\right)\right) =Re(\lambda)|\lambda|^2 \|w\|^2+\left(\alpha Re(\overline{\lambda}e^{-\lambda \tau})+|\lambda|^2 \beta \right)  \|w_{x}\|^2 + Re(\lambda)\|\varphi\|^2 + \kappa  \|\varphi_{
x}\|^2. $$
To conclude, it suffices to prove that $\left(\alpha Re(\overline{\lambda}e^{-\lambda \tau})+|\lambda|^2 \beta \right)>0$.

Notice that,
\begin{equation*}
\alpha Re(\overline{\lambda}e^{-\lambda \tau})+|\lambda|^2 \beta=\alpha e^{-a\tau}\left( a\cos(b\tau)-b\sin(b\tau)\right)+|\lambda|^2 \beta.\label{eqq}
\end{equation*}  
Let $s \in \mathbb{R}$ such that $\cos(s)=\frac{a}{|\lambda|}$ and $\sin(s)=\frac{b}{|\lambda|}$ we have
$$\alpha Re(\overline{\lambda}e^{-\lambda \tau})+|\lambda|^2 \beta=|\lambda|\left( \alpha e^{-a\tau}\cos(b\tau+s)+|\lambda| \beta\right) \geq |\lambda|\left( |\lambda| \beta-\alpha e^{-a\tau}\right)> |\lambda|\left( |\lambda| \beta-\alpha \right).$$
Then $\alpha Re(\overline{\lambda}e^{-\lambda \tau})+|\lambda|^2 \beta >0$ for $|\lambda|\geq \frac{\alpha}{\beta}$.

Now, we suppose that $|\lambda|\leq \frac{\alpha}{\beta}$.
We have
  \begin{eqnarray}
\alpha Re(\overline{\lambda}e^{-\lambda \tau})+|\lambda|^2 \beta &\geq & \alpha e^{-a\tau}\left( a\cos(b\tau)-b^2\tau\right)+(a^2+b^2)  \beta \notag
\\ 
&\geq & \alpha e^{-a\tau} a\cos(b\tau)+\left(\beta-\alpha\tau e^{-a\tau}\right)b^2+a^2\beta   . \label{COND} 
\end{eqnarray}   
  
By using $\alpha\tau \leq \beta$, we have $\left(\beta-\alpha\tau e^{-a\tau}\right)>0$, moreover, one has $\cos (b\tau) >0$ because  $|b\tau|\leq \frac{\tau \alpha}{\beta}\leq 1<\frac{\pi}{2}$. Then we deduce from (\ref{COND}) that
  $$\alpha Re(\overline{\lambda}e^{-\lambda \tau})+|\lambda|^2 \beta>0.$$

In conclusion, the sesquilinear form $B$ is coercive for every $\lambda \in \mathbb{C}_0$.

 By the Lax-Milgram lemma, equation (\ref{max99}) has a unique solution $(u,\theta ) \in \mathcal{F}.$ Immediately, from (\ref{max1sp1}$_1$), we have that $v \in H^1_0(\Omega).$ Now, if we consider $(w,\varphi) \in \lbrace 0\rbrace \times \mathcal{D}(\Omega)$ in (\ref{max99}) we deduce that  equation (\ref{max1sp1})$_4$ holds true and $\theta\in H^2(\Omega)$. The function $z$, defined by (\ref{max33}), belongs to $L^2(\Omega , H^1(0,1))$ and satisfies (\ref{max1sp1})$_3$ and $z(.,0)=u_x.$ 
 
The functions $z(.,1)$ and $v_x$ belong to $L^2(\Omega),$ then we take $(w,\varphi) \in \mathcal{D}(\Omega) \times \lbrace 0\rbrace$ in (\ref{max99}) to deduce that $\alpha z(.,1)+\beta v_x$ belongs to $H^1(\Omega)$ and that  equation (\ref{max1sp1})$_2$ holds true. 
 
 Let $\tilde{\theta} = \theta-\frac{1}{\ell}\int_\Omega \theta_0dx,$ then  we have that  $U=(u,v,z,\tilde{\theta})$ belongs to $\mathcal{D}(\mathcal{A})$, and $\mathcal{A}U=F.$ Thus, $\lambda I-\mathcal{A}$ is surjective
 
 To prove that $\lambda I-\mathcal{A}$ is injective, we solve the equation
 \begin{equation*}
(\lambda I-\mathcal{A})U=0, \label{sp1111}
\end{equation*}
where $U=(u,v,z,\theta) \in \mathcal{D}(\mathcal{A})$. Using the previous calculus (with $F=(0,0,0,0)$), we deduce easily that
$$Re(\lambda)|\lambda|^2 \|u\|^2+\left(\alpha Re(\overline{\lambda}e^{-\lambda \tau})+|\lambda|^2 \beta \right)  \|u_{x}\|^2 + Re(\lambda)\|\theta\|^2 + \kappa  \|\theta_{
x}\|^2=0$$
which implies that $U=0$. Thus $\lambda I-\mathcal{A}$ is injective. Finally, recall that $\mathcal{A}$ is closed, then $\lambda I-\mathcal{A}$ is boundedly invertible. Hence $\lambda \in \rho(\mathcal{A})$ (for every $\lambda\in\mathbb{C}_0$), where $\rho(\mathcal{A})$ is the resolvent set of $\mathcal{A}$.
\end{proof}
\begin{lemma} \label{cc2}
For $\xi=2\frac{\tau \alpha^2}{\beta}$  and $\alpha\tau\leq\beta$, it holds
\begin{equation}
\sup_{Re(\lambda)>0}\|(\lambda I-\mathcal{A})^{-1}\|_{\mathcal{L}(H)}<\infty, \label{cond2}
\end{equation}
\end{lemma}
\begin{proof}
Suppose that condition (\ref{cond2}) is false. Then there exists a sequence of complex numbers $\lambda_n$ such that $Re(\lambda_n)>0$, for all $n\in \mathbb{N}$, and $|\lambda_n|\rightarrow \infty$, and a sequence of vector $U_n=(u_n,v_n,z_n,\theta_n)\in\mathcal{D}(\mathcal{A})$ with $\|U_n\|=1$ for every $n\in\mathbb{N}$, such that
\begin{equation}
\|(\lambda_n I-\mathcal{A})U_n\|=o(1). \label{sp112}
\end{equation} 

i.e.
\begin{equation}
\left\{ 
\begin{tabular}{l}
$\lambda u_n- v_n=f_n=o(1),\;\;\text{in}\;\;H_0^1(\Omega),$ \\
$\lambda_n v_n-\left(\alpha z_n(.,1)+\beta v_{n,x} \right)_x+\gamma \theta_x =g_n=o(1),\;\;\text{in}\;\;L^2(\Omega) ,$\\
$\lambda_n z_n+ \frac{1}{\tau}z_{n,\rho} =p_n=o(1),\;\;\text{in}\;\;L^2(\Omega \times (0,1)),$ \\ 
$\lambda_n \theta_n + \gamma v_{n,x} -\kappa \theta_{n,xx}= h_n=o(1),\;\;\text{in}\;\;L^2(\Omega).$
\end{tabular}
\right.  \label{max1sp}
\end{equation}
We will prove that $\|U_n\|=o(1)$ which contradict the hypothesis: $\|U_n\|=1$.

Using (\ref{max1sp})$_1$, we deduce that
\begin{eqnarray*}
Re \left( \left\langle (\lambda_n I-\mathcal{ A})U_n,U_n\right\rangle _{\mathcal{H}}\right)&\geq & Re(\lambda_n)+ \frac{\beta}{2} \|v_{n,x}\|^2-m \|u_{n,x}\|^2+\kappa \|\theta _{n,x}\|^2\notag\\
&\geq & Re(\lambda_n)+ \frac{\beta}{2} \|v_{n,x}\|^2-\frac{m}{|\lambda|^2} \|v_{n,x}+f_{n,x}\|^2+\kappa \|\theta _{n,x}\|^2
\end{eqnarray*}
Using that $\|v_{n,x}+f_{n,x}\|^2 \leq 2\|v_{n,x}\|^2+2 \|f_{n,x}\|^2$, we get
\begin{equation}
Re \left( \left\langle (\lambda_n I-\mathcal{ A})U_n,U_n\right\rangle _{\mathcal{H}}\right)\geq  Re(\lambda_n)+ \left(\frac{1}{2}\beta-\frac{2m}{|\lambda_n|^2} \right)\|v_{n,x}\|^2+\kappa \|\theta _{n,x}\|^2-\frac{2m}{|\lambda_n|^2} \|f_{n,x}\|^2 \label{sss1}
\end{equation}
Since $|\lambda_n|\rightarrow\infty$, there exists $n_0 \geq 0$, such that for every $n\geq n_0$, $\frac{\beta}{2}-\frac{2m}{|\lambda_n|^2}\geq \frac{\beta}{4}>0$;
then using that $f_{n,x}=o(1)$, the estimates (\ref{sp112}) and (\ref{sss1}) imply that 
$$v_{n,x}=o(1)\;\;\;\text{and}\;\;\;\theta_{n,x}=o(1)$$
 which yield, first, using again (\ref{max1sp})$_1$,
$$u_{n,x}=o(1),$$  second, using Dirichlet conditions on $v_n$, and Neumann conditions on $\theta_n$ (taking in mind $\int_\Omega \theta_n dx=0$), we deduce that
$$v_{n}=o(1)\;\;\;\text{and}\;\;\;\theta_{n}=o(1).$$ 

To determine $z_n,$ recall that $z_n(.,0)=u_{nx},$ then, by $(\ref{max1sp})_3$, we obtain
\begin{equation*}
z_n(.,\rho)=e^{-\lambda_n \tau \rho}u_{n,x}+ \tau e^{-\lambda_n \tau \rho}\int_0^\rho p_n(s)e^{\lambda_n \tau s}ds,
\end{equation*}

then
$$\int_0^1 \|z_n(.,\rho)\|^2d\rho \leq 2\|u_{n,x}\|^2+2\tau^2\int_0^1 \|p_n(\rho)\|^2=o(1).$$
In conclusion, we proved that $\|U_n\|=o(1)$ which contradict the fact that $\|U_n\|=1$.
\end{proof}
In view of Lemma \ref{cc1} and Lemma \ref{cc2} the two conditions in Lemma  \ref{cc3} are satisfied by the operator $\mathcal{A}$, then the associated semigroup $e^{t\mathcal{A}}$ is exponentially stable. The proof of Theorem \ref{cc4} is then finished.
\end{proof}

\section*{Comments}
 We can replace the Neumann conditions for $\theta$
$$\theta_x(0,t)=\theta_x(\ell,t)=0$$ by the Dirichlet conditions
 $$\theta(0,t)=\theta(\ell,t)=0,$$ we then obtain the same results.
\section*{Statements and Declarations}
\textbf{Availability of data.} Not applicable.

\textbf{Conflict of interest.}
We declare that we have no conflict of interest.


\begin{thebibliography}{1}

\bibitem{Amm15'}
~Ammari K.  and Mercier D.: \emph{Boundary feedback stabilization of a chain of
  serially connected strings}. Evolution Equations and Control Theory. 
  \textbf{1}, 1-19 (2015)

\bibitem{KhSh21}
 Khatir S.~M. and Shel F.: \emph{Well-posedness and exponential stability of a
  thermoelastic system with internal delay}. Applicable Analysis. \textbf{101}, 4851-4865 (2022).
  Doi: 10.1080/00036811.2021.1873299

\bibitem{Liu99}
Liu Z. and Zheng S.: \emph{Semigroups associated with dissipative systems}.
  Chapman $\&$ Hall/CRC (1999)


\end{thebibliography}

\end{document}